\documentclass[11pt]{article}

\usepackage{amssymb,amsmath,amsthm,rotating}

\newcommand{\n}{\noindent}
\newcommand{\bb}[1]{\mathbb{#1}}
\newcommand{\cl}[1]{\mathcal{#1}}
\newcommand{\mf}[1]{\mathfrak{#1}}

\newcommand{\ovl}{\overline}

\theoremstyle{plain}

\newtheorem{lm}{Lemma}

\newtheorem{thm}{Theorem}
\newtheorem{cor}[thm]{Corollary}

\theoremstyle{definition}

\begin{document}

\title{On quasi-free Hilbert modules\thanks{
The research of both authors was supported in part by a DST-NSF, S\&T 
Cooperation Programme grant.\newline
The research was begun in July 2003 during a visit by the first author to IHES, 
funded by development leave from Texas A\&M University, and a visit by 
the second author to Paris VI, funded by a grant from IFCPAR. We thank both 
institutions for 
their hospitality.}}

\author{Ronald G.~Douglas and Gadadhar Misra}

\date{}
\maketitle

\begin{abstract}
In this note we settle some technical questions concerning finite rank 
quasi-free Hilbert modules and develop some useful machinery. In particular, we 
provide a method for determining when two such modules are unitarily equivalent. 
Along the way we obtain representations for module maps and study how to 
determine the underlying holomorphic structure on such modules.
\end{abstract}

\setcounter{section}{-1}

\section{Introduction}\label{sec0}

\indent

One approach to multivariate operator theory is via the study of Hilbert 
modules, 
which are Hilbert spaces that are acted upon by a natural algebra of 
functions holomorphic on some bounded domain in complex $n$-space ${\bb C}^n$, 
(cf.\ 
\cite{D-P}, \cite{C-G}). 
In this 
setting, concepts and techniques from commutative algebra as well as from 
algebraic and complex geometry can be used. In particular, general Hilbert  
modules can be studied using resolutions by simpler or more basic Hilbert 
modules. Such an approach generalizes the dilation theory studied in the one 
variable or single operator setting (cf.\ \cite{D-P}). In \cite{D-M} the 
existence of 
resolutions for a 
large class of Hilbert modules was established with the class of quasi-free 
Hilbert modules forming the building blocks. Such modules are defined as the 
Hilbert 
space completion of a space of vector-valued holomorphic functions that 
possesses a 
kernel 
function. It then follows that  a natural Hermitian holomorphic bundle is 
determined by such a module. However, for a given algebra there are many 
distinct, inequivalent Hilbert space completions, 
which raises the question of determining the relation between two such modules. 

In this note, we consider this question by examining more carefully the bundle 
associated with a quasi-free 
module and introduce a  non-negative matrix-valued modulus function for 
any pair of finite rank quasi-free Hilbert modules. We show that a necessary 
condition for the modules to be unitarily equivalent is for the modulus to be 
the absolute value of a holomorphic matrix-valued function. Moreover, if the 
domain is starlike, we show that this condition is also sufficient. The 
Hermitian 
holomorphic  vector bundle over $\Omega$ associated with a quasi-free Hilbert 
module possesses a 
natural 
connection and curvature. To prove our results we rely upon the localization 
characterization of unitary equivalence obtained in \cite{D-P}. In the rank one 
case, we have  line bundles and we show that the difference of the two 
curvatures is equal to 
the complex two-form-valued Laplacian of 
the logarithm of the modulus function. This identity enables one to 
reduce the question of unitary equivalence of two rank one quasi-free Hilbert 
modules to showing that the latter function vanishes identically.

Along the way we  examine closely how one obtains the holomorphic structure 
on the vector bundle defined by a quasi-free Hilbert module. To accomplish this 
we introduce the notion of kernel functions dual to a generating set and study 
concrete representations for module maps between two quasi-free Hilbert 
modules. These dual kernel functions are closely related to the usual 
two-variable kernel function. We also raise some related questions for more 
general Hilbert modules.

In our earlier work, we have assumed the algebra of functions is complete in the 
 supremum norm and hence that it is 
a commutative Banach algebra. While we continue to make that assumption in this 
note, we 
will 
point out along the way that much weaker assumptions are sufficient for many of 
the results. In 
particular,  when the domain is the unit ball, it is enough for the 
polynomial algebra to act on the Hilbert space so that the coordinate functions 
define
contraction operators.\medskip

\n {\bf Acknowledgment.} We want to thank Harold Boas and Mihai Putinar for some 
useful comments on the contents of this paper.

\section{The Modulus for Quasi-Free Hilbert Modules}\label{sec1}

\indent

We use kernel Hilbert spaces over bounded domains in 
${\bb C}^n$, which are also contractive Hilbert  modules for the natural 
function 
algebra over the domain. More precisely, we use the kind of Hilbert module 
introduced in 
\cite{D-M} for the study of module resolutions. We first recall the necessary 
terminology.

For $\Omega$ a bounded domain in $\bb C^n$, let $A(\Omega)$ be the function 
algebra obtained as the completion of the set of functions that are holomorphic 
in some neighborhood of the closure of $\Omega$. For $\Omega$ the unit ball $\bb 
B^n$ or the polydisk $\bb D^n$ in $\bb C^n$, we obtain the familiar ball and 
polydisk algebras, $A(\bb B^n)$ and $A(\bb D^n)$, respectively. The Hilbert 
space 
$\cl M$ is said to be a {\em contractive Hilbert module over ${\mit 
A(\Omega)}$\/} if 
$\cl M$ 
is a unital  module over $A(\Omega)$ with  module map $A(\Omega)\times \cl M\to 
\cl M$ such that
\[
\|\varphi f\|_{\cl M} \le \|\varphi\|_{A(\Omega)} \|f\|_{\cl M} \text{ for 
$\varphi$ in $A(\Omega)$ and $f$ in  } \cl M.
\]
The space $\cl R$ is said to be a {\em quasi-free Hilbert module of rank\/} $m$ 
{\em over\/} $A(\mit\Omega)$, $1\le 
m\le \infty$,  if it is obtained as the completion of the algebraic tensor 
product $A(\Omega) \otimes \ell^2_m$ relative to an inner product such that
\begin{itemize}
\item[1)] $\text{\em eval}_{\pmb{z}}\colon \ A(\Omega) \otimes \ell^2_m \to 
\ell^2_m$ is 
bounded for $\pmb{z}$ in $\Omega$ and locally uniformly bounded on $\Omega$;
\item[2)] $\|\varphi(\Sigma\theta_i\otimes x_i)\| = \|\Sigma\varphi\theta_i 
\otimes x_i\|_{\cl R} \le 
\|\varphi\|_{A(\Omega)} 
\| 
\Sigma 
\theta_i \otimes x_i\|_{\cl R}$ for $\varphi$, $\{\theta_i\}$ in $A(\Omega)$ and 
$\{x_i\}$ in $\ell^2_m$; and
\item[3)] for $\{F_i\}$ a sequence in $A(\Omega)\otimes \ell^2_m$ that is Cauchy 
in the 
$\cl 
R$-norm, it follows that $\text{\em eval}_{\pmb{z}}(F_i)\to 0$ for all $\pmb{z}$ 
in $\Omega$ 
iff 
$\|F_i\|_{\cl R}\to 0$.
\end{itemize}
\n Here, $\ell^2_m$ is the $m$-dimensional Hilbert space.

Actually, condition 2) can be replaced in this paper by:

\begin{itemize}
\item[2$'$)] $\|\varphi(\Sigma\theta_i\otimes x_i)\| \le 
K\|\varphi\|_{A(\Omega)} \|\Sigma\theta_i \otimes x_i\|_R$ for 
$\varphi, \{\theta_i\}$ in $A(\Omega)$ and $\{x_i\}$ in $\ell^2_m$ for some 
$K>0$.
\end{itemize}

Also, note that condition 3) already occurs in the fundamental paper of 
Aronszajn \cite{Arons} in which it is used to conclude that the abstract 
completion of a space of functions on some domain is again a space of functions.

There is another equivalent definition of quasi-free Hilbert module in terms of 
a generating set. The contractive Hilbert module ${\cl R}$ over $A(\Omega)$ is 
said to be 
quasi-free relative to the vectors $\{f_1,\ldots, f_m\}$ if the set generates 
${\cl R}$ and $\{f_i\otimes_A 1_z\}^m_{i=1}$ forms a basis for ${\cl R} 
\otimes_A {\bb C}_z$ for $\pmb{z}$ in $\Omega$. The set of vectors $\{f_i\}$ is 
called a {\em generating 
set\/} for ${\cl R}$. One must also assume that the evaluation functions 
obtained are locally uniformly bounded and that property 3) holds.
In \cite{D-M}, this characterization and other properties of quasi-free 
Hilbert 
modules are given. This concept is closely related to the notions of sharp and 
generalized Bergman kernels studied by Curto and Salinas \cite{C-S}, Agrawal  
 and Salinas 
 \cite{A-S}, and Salinas \cite{S}. We'll say more about this relationship later. 
Note that there is a significant difference between the notion of quasi-free and 
membership in class ${\cl B}_n(\Omega)$ in \cite{C-D2} and \cite{C-S}. For 
example, let ${\cl M}$ be the contractive Hilbert module over $A(\Delta)$ 
defined by the analytic Toeplitz operator $T_p$ on the Hardy space $H^2({\bb 
D})$, where the closure of $p({\bb D})$ is the closure of $\Delta$. Then ${\cl 
M}$ is in ${\cl B}_k(\Delta')$ for $\Delta'$ any domain in $\Delta$ disjoint 
from $p({\bb T})$, where $k$ is the winding number of the curve $p({\bb T})$ 
around $\Delta'$. However, ${\cl M}$ is a rank $k$ quasi-free Hilbert module 
relative to any algebra $A(\Delta')$ iff $p({\bb T})$ equals the boundary of 
$\Delta$, in which case $\Delta'=\Delta$ and $k$ is again the winding number.
 
We should mention that other authors have investigated the proper notion of 
freeness for topological modules over Frechet algebras (cf.\ pp.~76, 123 \cite{E-P}). 
Since one allows modules that are the direct sum of finitely many ? of the 
algebra or the topological tensor product of the algebra with a Frechet space, 
there can be a closer parallel with what is done in algebra.

Let $\cl R$ and $\cl R'$ each be a rank $m$ $(1\le m<\infty)$ quasi-free 
Hilbert module over $A(\Omega)$ for the generating sets of vectors $\{f_i\}$ and 
$\{g_i\}$, respectively. Then $\{f_i(\pmb{z})\}$ and $\{g_i(\pmb{z})\}$ each 
forms a basis 
for $\ell^2_m$ for $\pmb{z}$ on $\Omega$ and  $\cl R$ is the closure of the 
span of $\{\varphi f_i\mid \varphi\in A(\Omega), 1\le i\le m\}$ while $\cl R'$ 
is the closure of the span of $\{\varphi g_i\mid \varphi\in A(\Omega), 1\le i\le 
m\}$. Consider the subspace $\Delta$ of $\cl R\oplus\cl R'$ which is the closure 
of the linear span of $\{\varphi f_i\oplus \varphi g_i\mid \varphi\in A(\Omega), 
1\le 
i\le m\}$ in $\cl R\oplus\cl  R'$. Let $\text{\em Hol}_m(\Omega)$ be the space 
of 
all holomorphic $\mf L(\ell^2_m)$-valued functions on $\Omega$.

\begin{lm}\label{lm1}
The subspace $\Delta$ is the graph of a closed, densely defined, one-to-one 
transformation $\delta=\delta(\cl R,\cl R')$ having dense range. Moreover, the 
domain and range of 
$\delta$ are invariant under the module action and $\delta$ is a module 
transformation.
\end{lm}

\begin{proof}
Since $\Delta$ is closed and the domain and range of $\delta$, if it is 
well-defined, 
will contain the linear
spans of $\{\varphi f_i\mid \varphi \in A(\Omega), 1\le i\le m\}$ and $\{\varphi 
g_i\mid \varphi\in A(\Omega),1\le i\le m\}$, respectively, the only thing 
needing proof is 
that $h\oplus 0$ or $0\oplus k$ in $\Delta$ implies $h=0$ and $k=0$. For 
$0\oplus k$ in $\Delta$ we have sequences $\{\varphi^{(n)}_i\}$, $1\le i\le m$, 
such that $\Sigma \varphi^{(n)}_i f_i\to 0$, while $\Sigma \varphi^{(n)}_i 
g_i\to k$. Since evaluation at $\pmb{z}$ in $\Omega$ is continuous in the norm 
of $\cl R$, we have that $\Sigma\varphi^{(n)}_i(\pmb{z}) f_i(\pmb{z})\to 0$ for 
$\pmb{z}$ in $\Omega$. 
Since $\{f_i(\pmb{z})\}$ is a fixed basis for $\ell^2_m$, it follows that 
$\varphi^{(n)}_i(\pmb{z})\to 0$ for $1\le i\le m$. Hence, it follows that 
$k(\pmb{z})_ = \lim\limits_n\ \Sigma \varphi^{(n)}_i(\pmb{z}) g_i(\pmb{z}) = 0$ 
and since $k(\pmb{z}) = 0$ for $\pmb{z}$ in $\Omega$, we have $k=0$ by 3). The 
same 
argument works to show $h\oplus 0$ in $\Delta$ implies that $h=0$.
\end{proof}

Although the definition  of $\delta$ is given in terms of its graph for 
technical reasons, one should note that $\delta$ merely takes the given 
generating set 
for ${\cl R}$ to the generating set for ${\cl R}'$.

To consider the infinite rank case,  we would need to know more about the 
relationship between the sets of values of the generating sets 
$\{f_i(\pmb{z})\}$ and 
$\{g_i(\pmb{z})\}$ in 
$\ell^2_m$ for the preceding argument to succeed (cf.\ \cite{D-M}).

Note that the graph $\Delta$ can also be interpreted as a rank $m$ quasi-free 
Hilbert module over $A(\Omega)$ relative to the generating set $\{f_i\oplus 
g_i\}$. 
Moreover, if we repeat the above construction relative to the pairs 
$\{\Delta,\cl R\}$ and $\{\Delta, \cl R'\}$, the transformations 
$\delta(\Delta,\cl R)$ and 
$\delta(\Delta, \cl R')$ are bounded. Finally, since $\delta(\cl R,\cl R') = 
\delta(\Delta,\cl R')^{-1} \delta(\Delta,\cl R)$,  many calculations for 
$\delta(\cl R,\cl R')$ can be reduced to the analogous calculations for a 
bounded module map 
composed with the inverse of a bounded module map.

If evaluation on ${\cl R}$ and ${\cl R}'$ are both continuous, the lemma holds 
if 
we replace $A(\Omega)$ by any algebra of holomorphic functions $A$ so long as it 
is norm dense in $A(\Omega)$. For example, if $\Omega$ is the unit ball $\bb 
B^n$ or the polydisk $\bb D^n$, one could take $A$ to be the algebra of all 
polynomials $\bb C[\pmb{z}]$ or the algebra of functions holomorphic on some 
fixed neighborhood of the closure of $\Omega$.

Now recall that for $\pmb{z}$ in $\Omega$, one defines the module $\bb 
C_{\pmb{z}}$ over $A(\Omega)$, where $\bb C_{\pmb{z}}$ is the one-dimensional 
Hilbert space $\bb C$, such that $\varphi\times\lambda = 
\varphi(\pmb{z})\lambda$ for $\varphi$ in $A(\Omega)$ and $\lambda$ in $\bb 
C_{\pmb{z}}$. Note that $\cl R\otimes_{A(\Omega)} \bb 
C_{\pmb{z}}\cong\bb C_{\pmb{z}} \otimes\ell^2_m$ for $\cl R$ any rank $m$ 
quasi-free Hilbert module. 
Localization of a Hilbert module $\cl M$ at $\pmb{z}$ in $\Omega$ is defined to 
be the module tensor product $\cl M \otimes_{A(\Omega)} \bb C_{\pmb{z}}$ (cf.\ 
\cite{D-P}), which is canonically isomorphic to the quotient module $\cl M/\cl 
M_{\pmb{z}}$, where $\cl M_{\pmb{z}}$ is the closure of $A(\Omega)_{\pmb{z}}\cl 
M$ and $A(\Omega)_{\pmb{z}} = \{\varphi\in A(\Omega)\mid \varphi(\pmb{z}) = 
0\}$. 
(Again, we can define  this construction for an algebra $A$, as above, so long 
as 
the set of 
functions in $A$ that vanish at a fixed point $\pmb{z}$ in $\Omega$ is dense in 
$A(\Omega)_{\pmb{z}}$.)

In addition to localizing Hilbert modules, one can localize module maps. While 
localization of bounded module maps is 
straightforward, here we need to localize $\delta$  which is possibly unbounded 
and hence we must be somewhat 
careful.

\begin{lm}\label{lm2}
For $\pmb{z}$ in $\Omega$, the map $\delta\otimes_{A(\Omega)} 1_{\pmb{z}}\colon 
\ \cl R \otimes_{A(\Omega)} \bb C_{\pmb{z}} \longrightarrow \cl R' 
\otimes_{A(\Omega)} \bb C_{\pmb{z}}$ is well-defined. Moreover, $\delta 
\otimes_{A(\Omega)} 1_{\pmb{z}}$ is an invertible operator on the 
$m$-dimensional 
Hilbert space $\bb C_{\pmb{z}} \otimes \ell^2_m$.
\end{lm}

\begin{proof}
Since for $\pmb{z}$ in $\Omega$, $A(\Omega)_{\pmb{z}} f_i$ is contained in the 
domain of $\delta$ for $1\le i\le m$ and $\delta(A(\Omega)_{\pmb{z}} f_i)$ is 
contained in the linear span of
$\{A(\Omega)_{\pmb{z}} g_i\}$, $1\le i\le m$, we see that one can define 
$\delta$ 
from $\cl R/\cl 
R_{\pmb{z}}$ to $\cl R'/\cl R'_{\pmb{z}}$ as a densely defined, module 
transformation having dense range. Both
$\cl R/\cl R_{\pmb{z}}$ and $\cl R'/\cl R'_{\pmb{z}}$ are $m$-dimensional since 
they are isomorphic to $\cl R\otimes_{A(\Omega)}\bb C_{\pmb{z}}$ and $\cl 
R'\otimes_{A(\Omega)} \bb C_{\pmb{z}}$, respectively. Since $\delta$ has dense 
range, it follows that $\delta \otimes_{A(\Omega)} 1_{\pmb{z}}$ is onto and thus 
invertible. Therefore, the final statement 
holds.
\end{proof}

Localization as defined above is used implicitly in the work of Arveson and 
others. 
Consider, for example, the recent paper \cite{A} involving free covers. Since 
the defect space 
is 
simply $F \otimes_{{\bb C}[z]} {\bb C}_0$, the assumption in Definition 2.2 of  
\cite{A} is 
that the localization map $A \otimes_{{\bb C}[z]} I_z = \dot A$ is unitary. 
While this 
observation doesn't add anything per se, it does  raise the question about  
the meaning of localization at other $\pmb{z}$, not just at the origin. We'll 
say  more
about this matter later in this note. A similar 
question can be raised in the work of Davidson \cite{D} who uses the trace which 
is just the localization map from a module ${\cl M}$ to ${\cl M}\otimes_A {\bb 
C}_0$.  Does consideration of localization at other 
$\pmb{z}$ add anything? Since the algebra in this case is non-commutative, this 
question would likely take us into the realm of 
non-commutative algebraic geometry such as considered by Kontsevich and 
Rosenberg \cite{K-R}.

The modulus $\mu = \mu(\cl R,\cl R')$ of $\cl R$ and $\cl R'$ is defined to be 
the absolute value of $\delta\otimes_{A(\Omega)} 1_{\pmb{z}}$. For $m>1$, there 
are two possibilities:\ the 
square root of $(\delta \otimes_{A(\Omega)} 1_{\pmb{z}})^* (\delta 
\otimes_{A(\Omega)} 1_{\pmb{z}})$ and the square root of 
$(\delta\otimes_{A(\Omega)} 1_{\pmb{z}}) (\delta \otimes_{A(\Omega)} 
1_{\pmb{z}})^*$. The first operator, which we'll denote by $\mu(\cl R,\cl R')$, 
is defined on $\cl R\otimes_{A(\Omega)} \bb C_{\pmb{z}}$ while the second one, 
which corresponds to $\mu'({\cl R},{\cl R}')$, is defined on $\cl 
R'\otimes_{A(\Omega)} \bb C_{\pmb{z}}$. In either case, $\mu$ is an invertible 
positive 
$m\times m$ matrix function which is distinct from the absolute value of 
$\delta(\cl R',\cl R) = 
\delta(\cl R,\cl R')^{-1}$.

Next we need to know more about the adjoint transformation $\delta^*\colon \ \cl 
R'\to \cl R$. Recall we know from von~Neumann's fundamental results \cite{JN}, 
that 
$\delta^*$ exists and its graph is given by the orthogonal complement of 
$\Delta$, the graph of $\delta$, in $\cl R\oplus \cl R'$ after reversing the 
roles of $\cl R$ and $\cl R'$ and introducing a minus sign. In particular, the 
graph $\Delta^*$ of $\delta^*$ is equal to $\{h\oplus k \in \cl R'\oplus \cl 
R\mid -k\oplus h \perp \Delta\}$.

For $\pmb{z}$ in $\Omega$, let $\{k^i_{\pmb{z}}\}$ and $\{k'{}^i_{\pmb{z}}\}$ be 
elements in $\cl R$ and $\cl R'$, respectively, such that $\langle h(\pmb{z}), 
g_i(\pmb{z})\rangle_{\ell^2_m} = \langle h,k'{}^{i}_{\pmb{z}}\rangle_{\cl R'}$ 
and 
$\langle k(\pmb{z}),  f_i(\pmb{z})\rangle_{\ell^2_m} = \langle 
k,k^i_{\pmb{z}}\rangle_{\cl R}$ for $h$ and $k$ in $\cl R'$ and $\cl R$, 
respectively. Note that the sets $\{k^i_{\pmb{z}}\}$ and 
$\{k'{}^{i}_{\pmb{z}}\}$ 
span the orthogonal complements of $\cl R_{\pmb{z}}$ and $\cl R'_{\pmb{z}}$, 
respectively. We will refer to the sets $\{k^i_{\pmb{z}}\}$ and 
$\{k'{}^i_{\pmb{z}}\}$, as the {\em dual sets 
of kernel functions\/} for the generating sets $\{f_i\}$ for ${\cl R}$ and 
$\{g_i\}$ for 
${\cl R}'$, respectively. Finally, for 
$\pmb{z}$ in $\Omega$ let $X_{ij}(\pmb{z})$ be the 
matrix in $\mf L(\ell^2_m)$ that satisfies
\[
\left\langle\sum_j X_{ij}(\pmb{z}) f_j(\pmb{z}), 
f_\ell(\pmb{z})\right\rangle_{\ell^2_m} = \langle g_i(\pmb{z}), 
g_\ell(\pmb{z})\rangle_{\ell^2_m} \text{ for } 1\le i,\ell\le m.
\]
In other words, $\{X_{ij}\}$ effects the change of basis from $\{f_i\}$ for 
${\cl 
R}$ 
to $\{g_i\}$ for ${\cl R}'$.
If we define $Y(\pmb{z})\colon \ \ell^2_m\to \ell^2_m$ so that  
$Y(\pmb{z})f_i(\pmb{z}) = g_i(\pmb{z})$ for $1\le i\le m$, then $Y(\pmb{z})$ is 
invertible 
and $\{X_{ij}(\pmb{z})\}$ is the matrix defining the operator $Y(\pmb{z})^* 
Y(\pmb{z})$ on $\ell^2_m$.
Moreover, since the generating sets $\{f_i(\pmb{z})\}$ and $\{g_i(\pmb{z})\}$ 
are holomorphic, the matrix-function $X_{ij}(\pmb{z})$ is real-analytic.

\begin{lm}\label{lm3}
The domain of $\delta^*$ contains the finite linear span of 
$\{k'{}^{i}_{\pmb{z}}\mid \pmb{z}\in\Omega, 1\le i\le m\}$. Moreover,
\[
\delta^*k'{}^{i}_{\pmb{z}} = \sum_j X_{ij}(\pmb{z}) k^j_{\pmb{z}}.
\]
\end{lm}

\begin{proof}
Since the span of $\{\varphi f_i\oplus \varphi g_i\mid \varphi\in A(\Omega), 
1\le i\le m\}$ is dense in $\Delta$, it is enough to show that 
\[
\left\langle\left(-\sum_j X_{ij}(\pmb{z}) k^j_{\pmb{z}}\right) \oplus 
k'{}^{i}_{\pmb{z}}, 
\varphi f_\ell \oplus \varphi g_\ell\right\rangle = 0
\]
for $\varphi$ in $A(\Omega)$ and $1\le \ell \le m$. But
\begin{align*}
&\left\langle \left(-\sum_j X_{ij}(\pmb{z}) k^j_{\pmb{z}}\right)\oplus 
k'{}^{i}_{\pmb{z}}, 
\varphi f_\ell \oplus \varphi g_\ell\right\rangle_{\cl R\oplus\cl R'} = 
\left\langle -\sum_j X_{ij}(\pmb{z}) k^j_{\pmb{z}}, \varphi f_\ell 
\right\rangle_{\cl R} + \langle k'{}^{i}_{\pmb{z}}, \varphi g_\ell\rangle_{\cl 
R'}\\
&\qquad = -\sum_j X_{ij}(\pmb{z}) \overline{\varphi(\pmb{z})} \langle 
k^j_{\pmb{z}}, 
f_\ell\rangle_{\cl R} + \ovl{\varphi(\pmb{z})} \langle k'{}^{i}_{\pmb{z}}, 
g_\ell 
\rangle_{{\cl R}'}\\
&\qquad = \ovl{\varphi(\pmb{z})} \left(\left\langle-\sum_j X_{ij}(\pmb{z}) 
f_j(\pmb{z}), 
f_\ell(\pmb{z})\right\rangle_{\ell^2_m} + \langle g_i(\pmb{z}), 
g_\ell(\pmb{z})\rangle_{\ell^2_m} 
\right) = 0
\end{align*}
by the definition of $\{X_{ij}(\pmb{z})\}$ and thus
the result is proved.
\end{proof}

Before we proceed, the notion of the dual set of kernel functions can be used to 
establish the first notion of holomorphicity, or in fact in this case, 
anti-holomorphicity, of a quasi-free Hilbert module.

Suppose ${\cl R}$ is the completion of $A(\Omega) \otimes_{alg}\ell^2_m$ and we 
consider the generating set $\{1\otimes e_i\}$ for ${\cl R}$ with the dual set 
of kernel functions $\{k^i_{\pmb{z}}\}$. As we pointed out above, 
$\{k^i_{\pmb{z}}\}$ spans the orthonormal complement of ${\cl R}_{\pmb{z}}$ in 
${\cl R}$ for $\pmb{z}$ in $\Omega$. For $h$ in ${\cl R}$ we have $\langle 
k^i_{\pmb{z}},h\rangle_{\cl R} = \ovl{\langle h(\pmb{z}),e_i\rangle_{\ell^2_m}}$ 
which is an anti-holomorphic function on $\Omega$. Thus $k^i_{\pmb{z}}$ is a 
weakly anti-holomorphic function and therefore $\pmb{z}\longrightarrow 
k^i_{\pmb{z}}$ is strongly anti-holomorphic. Finally, since the functions 
$\{k^i_{\pmb{z}}\}$ span ${\cl R}^\bot_{\pmb{z}}$ for $\pmb{z}$ in $\Omega$, we 
see that $\bigcup\limits_{\pmb{z}\in\Omega} {\cl R}^\bot_{\pmb{z}}$ is an 
anti-holomorphic Hermitian rank $m$ vector bundle over $\Omega$.

We record this result as 

\begin{lm}\label{lm4}
For ${\cl R}$ a finite rank $m$ quasi-free Hilbert module, 
$\bigcup\limits_{\pmb{z}\in \Omega} {\cl R}^\bot_{\pmb{z}}$ is an 
anti-holomorphic Hermitian rank $m$ vector bundle over $\Omega$.
\end{lm}

With the additional assumption of a ``closedness of range'' condition, this 
result is established in \cite{C-S}. Also, the above proof can be rephrased in 
terms of the ordinary notion of kernel function and rests on the holomorphicity 
of the functions in ${\cl R}$. Note that we have assumed the local uniformed 
boundedness of evaluation to reach the conclusion of Lemma \ref{lm4}. It would 
be of interest to understand better the relation of this notion to that of the 
closedness of range condition. In particular, one knows that the latter property 
does not always hold although it is unclear  whether evaluation is always 
locally uniformly bounded.

\section{Representations of Module Maps}\label{sec2}

\indent

Next we state a result familiar in settings such as the one provided by that of 
quasi-free Hilbert modules, which we essentially used in the preceding section 
to define $\delta^*$.

\begin{lm}\label{lm5} 
If ${\cl R}$ and ${\cl R}'$ are quasi-free Hilbert modules over $A(\Omega)$ 
relative to the generating sets $\{f_i\}^m_{i=1}$ and $\{g_i\}^m_{i=1}$, $1\le m 
< 
\infty$, 
and $X$ is a module map from ${\cl R}$ to ${\cl R}'$, then there exists $\Psi= 
\{\psi_{ij}\}$ 
in  $\text{Hol}_m(\Omega)$ such that
\[
Xf_i = \sum^m_{j=1} \psi_{ij} g_j,\quad \text{for}\quad 1\le i\le m.
\]
\end{lm}

\begin{proof}
For $\pmb{z}$ in $\Omega$, both $\{f_i(\pmb{z})\}^m_{i=1}$ and 
$\{g_i(\pmb{z})\}^m_{i=1}$ are bases for $\ell^2_m$ and hence there exists a 
unique matrix $\{\psi_{ij}(\pmb{z})\}^m_{i,j=1}$ such that
\[
(Xf_i)(\pmb{z}) = \sum^m_{j=1} \psi_{ij}(\pmb{z}) g_j(\pmb{z})\quad 
\text{for}\quad i=1,2,\ldots, m.
\]
Since the functions $\{(Xf_i)(\pmb{z})\}^m_{i=1}$ and $\{g_i(\pmb{z})\}^m_{i=1}$ 
are all holomorphic, it follows from Cramer's rule that $\Psi = 
\{\psi_{ij}\}^m_{i,j=1}$ is in 
$\text{Hol}_m(\Omega)$ which completes the proof.
\end{proof}

Although we  obtain a holomorphic matrix function defining a module map 
between 
distinct quasi-free Hilbert 
modules,  this function is not very useful unless the 
modules 
and the generating sets are  
the same. That is because the matrix 
representing a linear transformation relative to different bases captures 
little information about the norm of it or the eigenvalues of its absolute 
value.

Before continuing, we want to show that the multiplier representation for a 
module map also extends to its localization.

\begin{lm}\label{lm6}
If ${\cl R}$ and ${\cl R}'$ are rank $m$ quasi-free Hilbert modules with 
generating sets 
$\{f_i\}$ and $\{g_i\}$, respectively, and $X\colon \ {\cl R}\to {\cl R}'$ is 
the 
module map  from ${\cl R}$ to ${\cl R}'$ represented by $\Psi = \{\psi_{ij}\}$ 
in $\text{Hol}_m(\Omega)$, then
\[
(X \otimes_A 1_{{\bb C}_{\pmb{z}}}) (f_i\otimes_A 1_{\pmb{z}}) = \sum^m_{j=1} 
\psi_{ij}(\pmb{z}) (g_j\otimes_A 1_{\pmb{z}})  \text{ for } \pmb{z} \text{ in } 
\Omega.
\]
\end{lm}

\begin{proof}
Let $\{k'{}^i_{\pmb{z}}\}$ be the  set of kernel functions dual to the 
generating set 
$\{g_i\}$. Then for a fixed $\pmb{z}$ the span of the set 
$\{k'{}^i_{\pmb{z}}\}^m_{i=1}$ is the 
orthogonal complement of $[A_{\pmb{z}}{\cl R}']$ and we can identify ${\cl R}' 
\otimes_A {\bb C}_{\pmb{z}}$ with the quotient module ${\cl R}'/[A_{\pmb{z}}{\cl 
R}']$. Calculating we see that the vector $Xf_i - \sum\limits^m_{i=1} 
\psi_{ji}(\pmb{z})g_j$ is orthogonal to each $k'{}^i_{\pmb{z}}$, $1\le \ell \le 
m$, and hence is in $[A_{\pmb{z}}{\cl R}']$. Therefore, we have that
\[
(X\otimes_A 1_{{\bb C}_{\pmb{z}}}) (f_i\otimes_A 1_{\pmb{z}}) = (Xf_i) \otimes_A 
1_{\pmb{z}} = \sum^m_{j=1} \psi_{ij}(z) (g_j\otimes_A 1_{\pmb{z}}) \quad 
\text{for}\quad 1\le i\le m,
\]
which completes the proof.
\end{proof}

Note that this result also holds for the localization of $\delta$. Also, if the 
ranks of ${\cl R}$ and ${\cl R}'$ are finite integers $m$ and $m'$ but not 
equal, then we 
obtain the same result for a holomorphic $m'\times m$ matrix-valued function.

Although, as we mentioned above, this representation has limited value, it does 
enable us to investigate the nature of the  sets of constancy for the local rank 
of a module map $X$ 
between two 
quasi-free Hilbert modules ${\cl R}$ and ${\cl R}'$. The previous lemma shows 
that, this local behavior is the same as that of a 
holomorphic matrix-valued function. In particular, the singular sets $\Sigma_k$ 
of $X \otimes_A 1_{\pmb{z}}$, that is, the subsets of $\Omega$ on which the 
rank of $X \otimes_A 1_{\pmb{z}}$ is $k$, are analytic subvarieties of $\Omega$. 
Thus we have established

\begin{thm}\label{thm1}
If ${\cl R}$ and ${\cl R}'$ are finite rank quasi-free Hilbert modules and $X$ 
is a module map $X\colon \ {\cl R}\to {\cl R}'$, then the singular sets 
$\Sigma_k$ of $X 
\otimes_A 1_{\pmb{z}}$ are analytic subvarieties of $\Omega$.
\end{thm}

We intend to use this fact to relate our work to that of Harvey--Lawson 
\cite{H-L} in the future. 
In particular, we expect their formulas for singular connections to be useful in 
obtaining invariants from resolutions such as those exhibited in \cite{D-M}.

\section{Holomorphic Structure}\label{sec3}

\indent

Recall that the spectral sheaf of a Hilbert module ${\cl M}$ over $A(\Omega)$ is 
defined to be $Sp({\cl M}) = \bigcup\limits_{\pmb{z}\in\Omega} {\cl M} \otimes_A 
{\bb C}_{\pmb{z}}$ with the collection of sections $\{f\otimes_A 1_{\pmb{z}}\mid 
f\in {\cl M}\}$. A priori the fibers of $Sp({\cl M})$ are isomorphic to the 
Hilbert modules ${\bb C}_{\pmb{z}} \otimes\ell^2_{{m}_{\pmb{z}}}$, where 
the dimension $m_{\pmb{z}}$ can vary from point to point and $0\le m_{\pmb{z}} 
\le \infty$. If ${\cl R}$ is a quasi-free rank $m$ Hilbert module, then 
$m_{\pmb{z}} = m$ for all $\pmb{z}$, but we would like more. Namely, we would 
like to define 
a canonical structure on $Sp({\cl R})$ making it into a holomorphic vector 
bundle relative to which the sections are holomorphic. We would also like to 
understand better the relation between the spectral sheaf $Sp({\cl R})$ and the 
anti-holomorphic vector bundle $\bigcup\limits_{\pmb{z} \in \Omega} {\cl 
R}^\bot_{\pmb{z}}$.

Although it might seem straightforward that the spectral sheaf $Sp({\cl R}) 
=\bigcup\limits_{\pmb{z}\in\Omega} {\cl R} \otimes_A {\bb C}_{\pmb{z}}$, for a 
finite rank quasi-free Hilbert module ${\cl R}$, is a Hermitian holomorphic 
vector bundle, it is worth considering how one exhibits such structure and shows 
that it is well-defined. 

Let $\{f_i\}^n_{i=1}$ be a subset of ${\cl R}$ relative to which ${\cl R}$ is 
quasi-free and define the map $F(\pmb{z})$ from ${\cl R} \otimes_A {\bb 
C}_{\pmb{z}}$ to $\ell^2_m$ such that $F(\pmb{z}) \left(\sum\limits^n_{i=1} 
\lambda_i(f_i \otimes_A 1_{\pmb{z}})\right) = \sum\limits^n_{i=1} \lambda_i 
f_i(\pmb{z})$. By the quasi-freeness of ${\cl R}$ relative to the generating set 
$\{f_i\}^m_{i=1}$, 
it follows that this map is well-defined, one-to-one and onto. Its inverse 
$F^{-1}$ 
defines a map from  the trivial vector bundle $\Omega \times \ell^2_m$ to the 
spectral sheaf $Sp({\cl R})$ 
of ${\cl R}$ which can be used to make $Sp({\cl R})$ into a holomorphic vector 
bundle. It is clear that the sections $f_i\otimes_{A(\Omega)} 1_{\pmb{z}}$ are 
holomorphic relative to this structure. We see later that the same is true for 
all $k$ in ${\cl R}$. The only issue now is whether the intrinsic norm on the 
fibers of 
$Sp({\cl R})$ yields a real-analytic metric on this bundle, which is necessary 
for 
$Sp({\cl R})$ to be a Hermitian holomorphic vector bundle.

To show that, consider $F(z)^{-1}\colon \ \ell^2_m\to {\cl R}\otimes_A {\bb 
C}_{\pmb{z}}$. We need to know that the function $\pmb{z}\to \langle F(z)^{-1}x, 
F(z)^{-1}y\rangle_{{\cl R} \otimes_A {\bb C}_{\pmb{z}}}$ is real-analytic for 
vectors $x$ and $y$ in $\ell^2_m$. Since the functions $\{f_i(\pmb{z})\}$ are 
holomorphic, the map from a fixed basis $\{e_i\}$ in $\ell^2_m$ to $\ell^2_m$ 
defined by $e_i \to f_i(\pmb{z})$ is holomorphic. Hence, the question rests on 
the behavior of the Grammian $\{\langle f_i\otimes_A 1_{\pmb{z}}, f_j \otimes_A 
1_{\pmb{z}}\rangle_{{\cl R} \otimes_A {\bf C}_{\pmb{z}}}\}$. Using the dual set 
of kernel functions $\{k^\ell_{\pmb{z}}\}^m_{\ell=1}$ for the generating set 
$\{f_i\}$, 
we see that $f_i \otimes_A 1_{\pmb{z}}$, viewed as a vector in ${\cl R}$, is the 
projection of $f_i$ onto ${\cl R}^\bot_{\pmb{z}}$, the span of the 
$\{k^\ell_{\pmb{z}}\}^m_{\ell=1}$. Now 
consider the identity involving the inner products $\langle 
f_i,k^\ell_{\pmb{z}}\rangle_{\cl 
R} = 
\langle f_i(\pmb{z}), f_\ell(\pmb{z})\rangle_{\ell^2_m}$ obtained using the 
defining 
property of the dual set $\{k^\ell_{\pmb{z}}\}$. We see that $\pmb{z}\to \langle 
f_i, k^\ell_{\pmb{z}}\rangle_{\cl R}$ is real-analytic. Therefore, inner 
products of the projections of $f_i$ and $f_j$ onto the span of the 
$\{k^\ell_{\pmb{z}}\}^m_{i=1}$ are also real-analytic which completes the proof. 
(Because of linear independence, the expressions can't vanish.)

 Now we must 
consider what happens if we use a different generating set $\{g_i\}^n_{i=1}$ 
relative to 
which ${\cl R}$ is quasi-free. Using  Lemma \ref{lm5}, we see that 
the map which sends $f_i$ to $g_i$, $i=1,2,\ldots,m$, is defined by a 
holomorphic $m\times 
m$ matrix-valued
 function $\Psi(\pmb{z})$ in $\text{Hol}_m(\Omega)$. That is, we have 
$g_i(\pmb{z}) = \sum\limits^m_{j=1} \psi_{ij}(\pmb{z}) f_j(\pmb{z})$ for 
$\pmb{z}$ in $\Omega$ and hence $\Psi(\pmb{z})$ defines a holomorphic bundle map 
which intertwines the holomorphic structures defined by the generating sets 
$\{f_i\}^n_{i=1}$ and 
$\{g_i\}^n_{i=1}$. Thus, we have proved:

\begin{thm}\label{thm2}
For ${\cl R}$ a finite rank quasi-free Hilbert module over $A(\Omega)$, there 
is a unique, well-defined holomorphic structure on $Sp({\cl R})$ relative to 
which the 
functions $\pmb{z}\to k \otimes_A 1_{\pmb{z}}$ are  holomorphic 
sections for each $k$ in ${\cl R}$.
\end{thm}

\begin{proof}
The only part requiring proof is the last statement. Clearly, this is true for 
any $f_i$ in a generating set $\{f_i\}^m_{i=1}$ for ${\cl R}$. Similarly, it 
follows 
for any  linear combination $\sum\limits^m_{i=1} \varphi_if_i$ for 
$\{\varphi_i\} \subset A(\Omega)$, that we obtain a holomorphic section. 
Finally, 
the ${\cl R}$-norm limit of such a sequence will converge uniformly locally and 
hence to a holomorphic section of $Sp({\cl R})$ which completes the proof.
\end{proof}

There is another approach to the holomorphic structure on $Sp({\cl R})$ which 
was essentially used in \cite{C-D2}, \cite{C-S}.
If the space $A_{\pmb{z}}{\cl R}$ is closed and the rank of ${\cl R}$ is finite, 
then the projection onto 
$[A_{\pmb{z}}{\cl R}]^\bot$ can be shown to define an anti-holomorphic map and 
hence the 
quotient ${\cl 
R}/[A_{\pmb{z}}{\cl R}]$ is holomorphic. Since ${\cl R}/[A_{\pmb{z}}{\cl R}] 
\cong {\cl R} \otimes_A {\bb C}_{\pmb{z}}$, this is another way of establishing 
a 
holomorphic structure on $Sp({\cl R})$. The smoothness of sections is 
straightforward in this case. However, the proof of Theorem \ref{thm2} is valid 
without the assumption of ``closed range''  but does require the local uniform 
boundedness of evaluation.

This identification of a holomorphic structure on the spectral sheaf of a finite 
rank quasi-free Hilbert module raises a series of questions regarding the 
situation for the spectral sheaf of a general Hilbert module. In particular, 
although we have called $Sp({\cl M}) = \bigcup\limits_{\pmb{z}\in\Omega} {\cl M} 
\otimes_A {\bb C}_{\pmb{z}}$ a sheaf, is it?

Although we can adopt the preceding approach to attempt to identify 
$\bigcup\limits_{\pmb{z}\in\Gamma} {\cl M}\otimes_A {\bb C}_{\pmb{z}}$ with 
the trivial bundle $\Gamma\times {\bb C}^m$ in case the fiber dimension is 
constant on an open subset 
$\Gamma$ of $\Omega$, the utility of this identification depends on being able 
to show that the transition functions on an overlap $\Gamma_1\cap \Gamma_2$ 
are holomorphic. This would show that $Sp({\cl M})$ is a holomorphic bundle for 
the 
``easy case,'' that is, a Hilbert module ${\cl M}$ for which the fiber dimension 
of ${\cl M}\otimes_A {\bb C}_{\pmb{z}}$ is constant and 
finite on all of $\Omega$. Until that case is decided, it is pointless to 
speculate about the general case of an ${\cl M}$ with finite but different 
dimensional fibers.

There is additional information about the behavior of the Grammian for the 
$\{f_i\otimes_A 1_{\pmb{z}}\}$ that we can obtain from a modification of the 
preceding arguments. Let $\{f_i\}$ be a generating set for the finite rank 
quasi-free Hilbert module ${\cl R}$.
We introduce a related notion of dual generating set which we will denote by 
$\{g^i_{\pmb{z}}\}$ so that $\langle h,g^i_{\pmb{z}}\rangle_{\cl R} = \langle 
h\otimes_A 1_{\pmb{z}}, f_i\otimes_A 1_{\pmb{z}}\rangle_{{\cl R} \otimes_A {\bb 
C}_{\pmb{z}}}$ for all $i$ and $\pmb{z}$ in $\Omega$ and $h$ in ${\cl R}$. If 
$P_{\pmb{z}}$ denotes the orthogonal projection of ${\cl R}$ onto ${\cl 
R}^\bot_{\pmb{z}}$, then one sees that $g^i_{\pmb{z}} = P_{\pmb{z}}f_i$ for all 
$i$ and $\pmb{z}$ in $\Omega$ since we can identify $f_i\otimes_A 1_{\pmb{z}}$ 
with $P_{\pmb{z}}f_i$. Since $\bigcup\limits_{\pmb{z}\in\Omega} {\cl 
R}^\bot_{\pmb{z}}$ is an anti-holomorphic Hermitian rank $m$ vector bundle, we 
see that the $\{g^i_{\pmb{z}}\}$ form an anti-holomorphic frame for it. 
Moreover, we have
\[
\langle f_i\otimes_A 1_{\pmb{z}}, f_j\otimes_A 1_{\pmb{z}}\rangle_{{\cl 
R}\otimes_A {\bb C}_{\pmb{z}}} = \langle P_{\pmb{z}} f_i, P_{\pmb{z}} 
f_j\rangle_{\cl R} = \langle g^i_{\pmb{z}}, g^j_{\pmb{z}}\rangle_{\cl R}
\]
or that the Grammian for the localization at $\pmb{z}$ in $\Omega$ of the 
generating set $\{f_i\}$ agrees with that of the anti-holomorphic frame 
$\{g^i_{\pmb{z}}\}$ for the anti-holomorphic Hermitian rank $m$ vector bundle 
$\bigcup\limits_{\pmb{z}\in\Omega} {\cl R}^\bot_{\pmb{z}}$. This allows us to 
obtain the following result which will be used in the next section.

\begin{thm}\label{thm3}
If ${\cl R}$ and ${\cl R}'$ are finite rank quasi-free Hilbert modules for the 
generating sets $\{f_i\}$ and $\{f'_i\}$ so that the Grammians $\{\langle f_i 
\otimes_A 1_{\pmb{z}}, f_j\otimes_A 1_{\pmb{z}}\rangle_{{\cl R}\otimes_A {\bb 
C}_{\pmb{z}}}\}$ and $\{\langle f'_i\otimes_A 1_{\pmb{z}}, f'_j \otimes_A 
1_{\pmb{z}}\rangle_{{\cl R}' \otimes_A {\bb C}_{\pmb{z}}}\}$ are equal, then 
$\delta({\cl R}, {\cl R}')$ is an isometric module map and ${\cl R}$ and ${\cl 
R}'$ are unitary equivalent.
\end{thm}

\begin{proof}
Proceeding as above we obtain anti-holomorphic frames $\{g^i_{\pmb{z}}\}$ and 
$\{g'{}^i_{\pmb{z}}\}$ for $\bigcup\limits_{\pmb{z}\in\Omega}{\cl 
R}^\bot_{\pmb{z}}$ and $\bigcup\limits_{\pmb{z} \in\Omega} {\cl 
R}'^\bot_{\pmb{z}}$, respectively. The mapping taking one anti-holomorphic frame 
to the other 
defines an anti-holomorphic unitary bundle map, call it $\Psi$, and hence the 
bundles are 
equivalent. Appealing to the Rigidity Theorem in \cite{C-D2}, we obtain a 
unitary 
operator $U\colon \ {\cl R}\to {\cl R}'$ which agrees with the bundle map, that 
is, $\Psi(\pmb{z}) = P'_{\pmb{z}} U|_{{\cl R}^\bot_{\pmb{z}}}$ for $\pmb{z}$ in 
$\Omega$. Moreover, since the action of $M^*_\varphi$ on ${\cl 
R}^\bot_{\pmb{z}}$ and ${\cl R}'{}^\bot_{\pmb{z}}$ is multiplication by 
$\ovl{\varphi(\pmb{z})}$, where $M_\varphi$ denotes the module actions of 
$\varphi$ on ${\cl R}$ and ${\cl R}'$, respectively, we see that $U^*$ is a 
module map from ${\cl R}'$ to ${\cl R}$ and hence $U = (U^*)^{-1}$ is a module 
map, which concludes the proof.
\end{proof}

\section{Equivalence of Quasi-Free Hilbert Modules}\label{sec4}

\indent

We now state our  first result about equivalence and the modulus..

\begin{thm}\label{thm4}
If the finite rank quasi-free Hilbert modules ${\cl R}$ and ${\cl R}'$ over 
$A(\Omega)$ are unitarily equivalent, then the modulus $\mu({\cl R}, {\cl R}')$ 
is the absolute value of a function $\Psi$ in $\text{Hol}_m(\Omega)$.
\end{thm}

\begin{proof}
Let $V\colon \ {\cl R}'\to {\cl R}$ be a unitary module map. We consider 
localization  of the triangle
\[
\begin{matrix}
\cl R \otimes_{A(\Omega)}\bb C_{\pmb{z}}&\overset{(V\delta)\otimes_{A(\Omega)} 
1_{\pmb{z}}}{\hbox to 45pt{\rightarrowfill}}&\cl R \otimes_{A(\Omega)} \bb 
C_{\pmb{z}}\\
\delta\otimes_{A(\Omega)} 1_{\pmb{z}}
\raisebox{2ex}{$\begin{rotate}{-45}{\hbox to 
25pt{\rightarrowfill}}\end{rotate}$}\hfill&&\raisebox{-2ex}{$\begin{rotate}{45}{
\hbox to 25pt{\rightarrowfill}}\end{rotate}$}~~~~~
 V\otimes_{A(\Omega)} 1_{\pmb{z}}\hfill\\
&\cl R' \otimes_{A(\Omega)} {\bb C}_{\pmb{z}}
\end{matrix}
\]
which yields $(V\delta)\otimes_{A(\Omega)} 1_{\pmb{z}} = (V\otimes_{A(\Omega)} 
1_{\pmb{z}}) (\delta \otimes_{A(\Omega)} 1_{\pmb{z}})$. Since $(V\delta) 
\otimes_{A(\Omega)} 1_{\pmb{z}}$ is  in $\text{\em 
Hol}_m(\Omega)$ by Lemmas \ref{lm5} and \ref{lm6}, it is sufficient to show that 
$V 
\otimes_{A(\Omega)} 1_{\pmb{z}}$ is unitary.

Again, by considering the factorization $I_{\cl R} \otimes_{A(\Omega)} 
1_{\pmb{z}} = (V^{-1} \otimes_{A(\Omega)} 1_{\pmb{z}}) (V\otimes_{A(\Omega)} 
1_{\pmb{z}})$ and in view of the fact that both $\|V^{-1} \otimes_{A(\Omega)} 
1_{\pmb{z}}\|\le \|V^{-1}\| = 1$ and $\|V\otimes_{A(\Omega)} 1_{\pmb{z}}\| \le 
\|V\| = 1$, we see that $V \otimes_{A(\Omega)} 1_{\pmb{z}}$ is unitary and the 
result is proved since $\mu(\cl R,\cl R')$ is the absolute value of $\delta(\cl 
R,\cl R')$.
\end{proof}

Note that if we use $V^{-1}$ from ${\cl R}$ to ${\cl R}'$ we see that the other 
square 
root, $\mu(\cl R',\cl R)$ is also the modulus of a holomorphic function in 
$\text{\em Hol}_m(\Omega)$.

The argument in this theorem raises a question about  a bounded module 
map $V$ between finite rank, quasi-free Hilbert module ${\cl R}'$ and ${\cl R}$ 
such 
that 
the localization $V\otimes_{\rm A(\Omega)} 1_{\pmb{z}}$ is unitary for $\pmb{z}$ 
in $\Omega$. We see by Theorem \ref{thm3} that such a map must be unitary if it 
has dense range by choosing a generating set $\{f_i\}$ for ${\cl R}'$ and the 
generating set $\{Vf_i\}$ for ${\cl R}$. If $\theta$ is a singular inner 
function, then the module map from the Hardy module $H^2({\bb D})$ to itself 
defined by multiplication by $\theta$ is locally one to one but does not have 
dense range. However, it is not locally a unitary map. It would seem likely that 
 maps that are locally unitary must have dense range but we have been unable to 
prove this. Some of 
these issues 
would also seem to be related to the proof of Theorem 2.4 in \cite{A}. This is 
the reference  we made earlier to the use in this work of localization at 
$\pmb{z}$ in addition to the origin.

What about the converse to the theorem? Suppose there exists a function $\Psi$ 
in $\text{\em 
Hol}_m(\Omega)$ such that $\Psi(\pmb{z})^* \Psi(\pmb{z}) = \mu(\pmb{z})^2 = 
(\delta \otimes_{A(\Omega)} 1_{\pmb{z}})^* (\delta\otimes_{A(\Omega)} 
1_{\pmb{z}})$. Since $\mu(\pmb{z})$ is invertible, we see that 
$\Psi(\pmb{z})^{-1}$ 
exists. 
Multiplying on the left by $(\Psi(\pmb{z})^{-1})^*$ and on the right by 
$\Psi(\pmb{z})^{-1}$, we obtain
\[
I = [(\delta \otimes_{A(\Omega)} 1_{\pmb{z}}) \Psi(\pmb{z})^{-1}]^* = [(\delta 
\otimes_{A(\Omega)} 1_{\pmb{z}}) \Psi(\pmb{z})^{-1}].
\]
Thus the function $(\delta \otimes_{A(\Omega)} 1_{\pmb{z}}) \Psi(\pmb{z})^{-1} = 
U(\pmb{z})$ is unitary-valued. We would like to show under these circumstances 
that 
$\cl R$ 
and $\cl R'$ are unitarily equivalent. The obvious approach is to consider the 
operator on $\cl R$ defined to be multiplication by $\Psi^{-1}$ followed by 
$\delta$. Unfortunately, we know little about the growth of $\Psi^{-1}$ as a 
function of $\pmb{z}$ and 
hence we don't know if the operator defined by multiplication by $\Psi$ is 
densely 
defined.

Suppose we assume that $\Omega$ is starlike relative to the point 
$\pmb{\omega}_0$ in $\Omega$, that is, the line segment $\{t\pmb{\omega}_0 + 
(1-t)\pmb{\omega}\mid 0\le t\le 1\}$ is contained in $\Omega$ for each 
$\pmb{\omega}$ in $\Omega$. Without loss of generality, we can assume that 
$\pmb{\omega}_0 = \pmb{0}$. Then we can define the  function $\Psi^{-1}_t \colon 
\ \Omega\to \mf L(\ell^2_m)$ for $0<t\le 1$ by $\Psi^{-1}_t(\pmb{z}) = 
\Psi^{-1}(t\pmb{z})$ for $\pmb{z}$ in $\Omega$. Now the family $\{\Psi^{-1}_t\}$ 
converge uniformly to $\Psi^{-1}$ on compact subsets of $\Omega$. (Actually, not 
only do the functions, which comprise the matrix entries, converge but so do all 
of their partial derivatives converge on compact subsets of $\Omega$.) 
Moreover, the matrix entries for $\{\Psi^{-1}_t\}$ for $0<t<1$ are in 
$A(\Omega)$ and thus we can define multiplication by $\Psi^{-1}_t$ on $\cl R$ 
and also $\delta\Psi^{-1}_t$. Moreover, $\delta\Psi^{-1}_t$ is a closed module 
transformation which has the same domain and range as $\delta$.

\begin{thm}\label{thm5}
If $\Omega$ is starlike and the modulus $\mu(\cl R,\cl R')$ for two finite 
rank quasi-free Hilbert modules over $A(\Omega)$ is the absolute value of a 
function in $\text{Hol}_m(\Omega)$, then $\cl R$ and $\cl R'$ are unitarily 
equivalent.
\end{thm}

\begin{proof}
By Lemma \ref{lm2} the localizations of both $\delta$ and $\delta\Psi^{-1}_t$ 
are 
well-defined and can be evaluated using the identifications of $\cl 
R\otimes_{A(\Omega)} \bb C_{\pmb{z}}$ and $\cl R'\otimes_{A(\Omega)} \bb 
C_{\pmb{z}}$ with $\cl R/\cl R_{\pmb{z}}$ and $\cl  R'/\cl R'_{\pmb{z}}$, 
respectively. For $\Phi$ a function in $\text{\em Hol}_m(\Omega)$ with entries 
from $A(\Omega)$, the operator $M_\Phi$ in $\mf L(\cl R)$ defined to be 
multiplication by $\Phi$, using generating sets for ${\cl R}$ and ${\cl R}'$, is 
well-defined and $M_\Phi \otimes_{A(\Omega)} 
1_{\pmb{z}} = \Phi(\pmb{z})$ for $\pmb{z}$ in $\Omega$. Next we consider the 
localization of the factorization of $\delta\Psi^{-1}_t$ to obtain
\begin{align*}
(\delta\Psi^{-1}_t) \otimes_{A(\Omega)} 1_{\pmb{z}} &= (\delta 
\otimes_{A(\Omega)} 1_{\pmb{z}}) (\Psi^{-1}_t \otimes_{A(\Omega)} 1_{\pmb{z}})\\
&= (\delta \otimes_{A(\Omega)} 1_{\pmb{z}}) \Psi^{-1}_t(\pmb{z})\\
&= U(\pmb{z}) [\Psi(\pmb{z}) \Psi^{-1}_t(\pmb{z})].
\end{align*}
Since $U(\pmb{z}) = (\delta \otimes_{A(\Omega)} 1_{\pmb{z}}) \Psi^{-1}(\pmb{z})$ 
is unitary, we see that the map $(\delta\Psi^{-1}_t)\otimes_{A(\Omega)} 
1_{\pmb{z}}$, which acts between the local modules $\cl R\otimes_{A(\Omega)} \bb 
C_{\pmb{z}}$ and $\cl R' \otimes_{A(\Omega)} \bb C_{\pmb{z}}$, is almost a 
unitary module map. Since $\lim\limits_{t\to 1} [\Psi(\pmb{z}) 
\Psi^{-1}_t(\pmb{z})] = I_{\ell^2_m}$, we see that the two local modules are 
unitarily 
equivalent. But for $m>1$ this is not enough. 

For $\cl M$ a Hilbert module and $n$ a positive integer, let $\cl 
M^n_{\pmb{z}}$ denote the closure of $(A(\Omega)^n_{\pmb{z}})\cl M$, where 
$A(\Omega)^n_{\pmb{z}}$ is the ideal of $A(\Omega)$ generated by the products 
of $n$ functions in $A(\Omega)_{\pmb{z}}$. (The quotient $\cl M/\cl 
M^n_{\pmb{z}}$ can also be identified as the module tensor product of $\cl M$ 
with some finite dimensional module with support at $\pmb{z}$. It is not 
straightforward, however, to identify the correct norm on the local module.) In 
Theorem 3.12 \cite{C-D}, X.\ Chen 
and 
the first author established that  a class of Hilbert modules, which 
includes the 
finite rank quasi-free Hilbert modules, are determined up to unitary equivalence 
by the collection of local modules $\cl M/\cl M^n_{\pmb{z}}$ for $\pmb{z}$ in 
$\Omega$, where $n$ depends on the 
rank of ${\cl R}$. To apply this result to $\cl R$ and $\cl R'$ we require the 
unitary equivalence of the higher order local modules $\cl R/\cl 
R^n_{\pmb{z}}$ and $\cl R'/\cl R'{}^n_{\pmb{z}}$. This is accomplished by noting 
that the 
localization of $[\Psi(\pmb{z}) \Psi^{-1}_t(\pmb{z})]$ to ${\cl R}'/{\cl 
R}'{}^n_{\pmb{z}}$ depends on the values of 
the partial derivatives of the entries of this matrix function up to some fixed 
order 
depending on $n$. Since the latter functions
all converge to the appropriate entries for the identity matrix on ${\cl 
R}'/{\cl R}'{}^n_{\pmb{z}}$, we conclude 
that $\cl 
R/\cl R^n_{\pmb{z}}$ and $\cl  R'/\cl R'{}^n_{\pmb{z}}$ are unitarily 
equivalent as $A(\Omega)$-modules. Thus, we conclude that $\cl R$ and $\cl R'$ 
are unitarily equivalent as $A(\Omega)$-modules.
\end{proof}

Actually $\Omega$ being starlike is not necessary. What is required for the 
preceding argument to work is that one 
can approximate the function $\Psi$ by matrix functions with entries from 
$A(\Omega)$ in a very strong sense. That is, one must be able to control not 
only the convergence of the function entries but also the 
convergence  of their partial derivatives and their inverses. By Montel's 
Theorem uniform convergence on compact subsets of $\Omega$ is sufficient. One 
can show using various techniques (cf.\ \cite{H} and Thm.\ 3.5.1 in \cite{HL}) 
that such approximation is 
possible for $\Omega$  bounded strongly pseudo-convex domain which allows us to 
state:

\begin{cor}\label{cor6}
If $\Omega$ is a bounded strongly pseudo-convex domain in ${\bb C}^m$ and the 
modulus $\mu(\cl R,\cl R')$ for two finite 
rank quasi-free Hilbert modules over $A(\Omega)$ is the absolute value of a 
function in $\text{Hol}_m(\Omega)$, then $\cl R$ and $\cl R'$ are unitarily 
equivalent.
\end{cor}

If we actually know that the mapping $\delta\Psi^{-1}$ is densely defined, we 
can use Theorem \ref{thm3} which means appealing to the Rigidity Theorem of 
\cite{C-D2} rather 
than involving curvature and its partial derivatives.

Now one knows that  a non-negative real-valued function $h(\pmb{z})$ on as 
simply connected domain $\Omega$ 
is 
the 
absolute value of a function holomorphic on $\Omega$ if and only if the 
two-form-valued  Laplacian of the logarithm of it vanishes identically on 
$\Omega$. 
Hence, we could restate Theorems \ref{thm4} and \ref{thm5} for the rank one case 
using this fact. However, we can go even 
further.

Recall we saw in Theorem \ref{thm2} that a rank $m$ quasi-free Hilbert module 
$\cl R$ 
determines a Hermitian 
holomorphic rank $m$ vector bundle $E_{\cl R} = 
\bigcup\limits_{\pmb{z}\in\Omega} \cl R\otimes_{A(\Omega)} \bb C_{\pmb{z}}$ 
over $\Omega$. Moreover, on such a bundle there is a canonical connection and 
hence a curvature which is a two-form valued matrix function on $\Omega$ (cf.\ 
\cite{C-D2}). In
the rank one case, we obtain a line bundle and if $\gamma(\pmb{z})$ is the 
holomorphic section $f\otimes_{A(\Omega)} 1_{\pmb{z}}$ of it, then the curvature 
$K_{\cl R}$ can be calculated so 
that
\[
 K_{\cl R}(\pmb{z}) = -\frac12 \sum_{i,j}\frac{\partial^2}{\partial 
z_i\partial \bar z_j} \log \|\gamma(\pmb{z})\| dz_i\wedge d\bar z_j.
\]

Now let us return to the case of two rank one quasi-free Hilbert modules over 
$\Omega$. If $\gamma'(\pmb{z})$ is the holomorphic section $g\otimes_{A(\Omega)} 
1_{\pmb{z}}$ for $E_{{\cl 
R}'}$, then $(\delta\gamma)(\pmb{z})$ is the 
holomorphic section $\gamma'(\pmb{z})$ for ${\cl 
R}'\otimes_{A(\Omega)} \bb C_{\pmb{z}}$. Moreover, a calculation shows that
\[
\|\gamma'(\pmb{z})\| = \|(\delta\gamma)(\pmb{z})\| = |(\delta 
\otimes_{A(\Omega)} 1_{\pmb{z}})| \|\gamma(\pmb{z})\|.
\]

\begin{thm}\label{thm6}
If $\cl R$ and $\cl R'$ are rank one quasi-free Hilbert modules and $\mu$ is the 
modulus, $\mu({\cl R},{\cl R}')$, then
\[
-\frac12 \sum_{i,j} \frac{\partial^2}{\partial z_i\partial \bar z_j} 
\mu(\pmb{z}) dz_i\wedge d\bar z_j =  
K_{\cl R} - K_{\cl R'}.
\]
\end{thm}

\begin{proof}
If $\gamma(\pmb{z})$ and $\gamma'(\pmb{z})$ are the holomorphic sections of 
$E_{\cl R}$ and $E_{{\cl 
R}'}$ given above, then  we have
\[
 K_{\cl R}  = -\frac12 \sum_{i,j} \frac{\partial^2}{\partial  z_i 
\partial \bar z_j} \log \|\gamma(\pmb{z})\| dz_i\wedge d\bar z_j \text{ 
and }  K_{\cl R'} 
= -\frac12 \sum_{i,j} \frac{\partial^2}{\partial  z_i \partial\bar z_j} \log 
|\delta\otimes_{A(\Omega)}1_{\pmb{z}}|
\|\gamma(\pmb{z})\| dz_i\wedge d\bar z_j.
\]
The proof is completed by using Lemma \ref{lm5} to conclude that $\mu(\pmb{z}) = 
|(\delta 
\otimes_{A(\Omega)} 1_{\pmb{z}})|$ for $\pmb{z}$ in $\Omega$.
\end{proof}

Formulas such as this one appeared first for specific examples in \cite{D-P} and 
for general quotient modules in \cite{D-M2} where they are used to obtain 
invariants for the quotient module. Here, of course, there is no quotient 
involved.

Finally, one can rephrase this result to state that for rank one quasi-free 
Hilbert modules the modulus is the square of the absolute value of a holomorphic 
function if and only if their respective curvatures coincide.

\end{document}